\documentclass[12pt,a4paper]{article}
 \usepackage[colorlinks]{hyperref}
\usepackage{a4,amsmath,amssymb,amsfonts,amsthm}
\usepackage{amsrefs}
\usepackage[all]{xy}
\usepackage[small,nohug,heads=vee]{diagrams}

\begin{document}


\newcommand{\fr}{\selectlanguage{francais}}
\newcommand{\en}{\selectlanguage{english}}


\def\M#1{\mathbb#1} 
\def\mR{\M{R}}           
\def\mZ{\M{Z}}           
\def\mN{\M{N}}           
\def\mQ{\M{Q}}       
\def\mC{\M{C}}  
\def\mG{\M{G}}
\def\mP{\M{P}}


\def\Spec{{\rm Spec}}
\def\rg{{\rm rg}}
\def\Hom{{\rm Hom}}
\def\Aut{{\rm Aut}}
 \def\Tr{{\rm Tr}}
 \def\Exp{{\rm Exp}}
 \def\Gal{{\rm Gal}}
 \def\End{{\rm End}}
 \def\det{{{\rm det}}}
 \def\Td{{\rm Td}}
 \def\ch{{\rm ch}}
 \def\che{{\rm ch}_{\rm eq}}
  \def\Spec{{\rm Spec}}
\def\Id{{\rm Id}}
\def\Zar{{\rm Zar}}
\def\Supp{{\rm Supp}}
\def\eq{{\rm eq}}
\def\Ann{{\rm Ann}}
\def\LT{{\rm LT}}
\def\Pic{{\rm Pic}}
\def\rg{{\rm rg}}
\def\et{{\rm et}}
\def\sep{{\rm sep}}
\def\ppcm{{\rm ppcm}}
\def\ord{{\rm ord}}
\def\Gr{{\rm Gr}}
\def\rk{{\rm rk}}
\def\Stab{{\rm Stab}}
\def\im{{\rm im}}
\def\Sm{{\rm Sm}}
\def\red{{\rm red}}
\def\Frob{{\rm Frob}}
\def\Ver{{\rm Ver}}
\def\Div{{\rm Div}}
\def\ab{{\rm ab}}


\def\beginProof{\par{\bf Proof. }}
 \def\endProof{${\qed}$\par\smallskip}
 \def\pr{^{\prime}}
 \def\prpr{^{\prime\prime}}
 \def\mtr#1{\overline{#1}}
 \def\ra{\rightarrow}
 \def\mfp{{\mathfrak p}}
 \def\mfm{{\mathfrak m}}
 
 \def\char{{\rm char}}
  \def\mQ{{\Bbb Q}}
 \def\mR{{\Bbb R}}
 \def\mZ{{\Bbb Z}}
 \def\mC{{\Bbb C}}
 \def\mN{{\Bbb N}}
 \def\mF{{\Bbb F}}
 \def\mA{{\Bbb A}}
  \def\mG{{\Bbb G}}
 \def\CI{{\cal I}}
 \def\CA{{\cal A}}
 \def\CD{{\cal D}}
 \def\CE{{\cal E}}
 \def\CJ{{\cal J}}
 \def\CH{{\cal H}}
 \def\CO{{\cal O}}
 \def\CA{{\cal A}}
 \def\CB{{\cal B}}
 \def\CC{{\cal C}}
 \def\CK{{\cal K}}
 \def\CL{{\cal L}}
 \def\CI{{\cal I}}
 \def\CM{{\cal M}}
  \def\CN{{\cal N}}
\def\CP{{\cal P}}
 \def\CZ{{\cal Z}}
\def\CR{{\cal R}}
\def\CG{{\cal G}}
\def\CX{{\cal X}}
\def\CY{{\cal Y}}
\def\CV{{\cal V}}
\def\CW{{\cal W}}
\def\CT{{\cal T}}
 \def\wt#1{\widetilde{#1}}
 \def\mod{{\rm mod\ }}
 \def\refeq#1{(\ref{#1})}
 \def\blb{{\big(}}
 \def\brb{{\big)}}
\def\mc{{{\mathfrak c}}}
\def\mcpr{{{\mathfrak c}'}}
\def\mcprpr{{{\mathfrak c}''}}
\def\ss{{\rm ss}}
\def\parf{{\rm parf}}
\def\P1{{{\bf P}^1}}
\def\cod{{\rm cod}}
\def\pr{\prime}
\def\prpr{\prime\prime}
\def\ss{\scriptstyle}
\def\OX{{ {\cal O}_X}}
\def\mpartial{{\mtr{\partial}}}
\def\inv{{\rm inv}}
\def\indlim{\underrightarrow{\lim}}
\def\prolim{\underleftarrow{\lim}}
\def\pprolim{'\prolim'}
\def\Pro{{\rm Pro}}
\def\Ind{{\rm Ind}}
\def\Ens{{\rm Ens}}
\def\without{\backslash}
\def\pbdb{{\Pro_b\ D^-_c}}
\def\qc{{\rm qc}}
\def\Com{{\rm Com}}
\def\an{{\rm an}}
\def\gfield{{\rm\bf k}}
\def\s{{\rm s}}
\def\dR{{\rm dR}}
\def\ari#1{\widehat{#1}}
\def\ul#1{\underline{#1}}
\def\sul#1{\underline{\scriptsize #1}}
\def\mou{{\mathfrak u}}
\def\ich{\mathfrak{ch}}
\def\cl{{\rm cl}}
\def\K{{\rm K}}
\def\R{{\rm R}}
\def\F{{\rm F}}
\def\L{{\rm L}}
\def\pgcd{{\rm pgcd}}
\def\rc{{\rm c}}
\def\N{{\rm N}}
\def\E{{\rm E}}
\def\H{{\rm H}}
\def\CHOW{{\rm CH}}
\def\A{{\rm A}}
\def\d{{\rm d}}
\def\Res{{\rm  Res}}
\def\GL{{\rm GL}}
\def\Alb{{\rm Alb}}
\def\alb{{\rm alb}}
\def\Hdg{{\rm Hdg}}
\def\Num{{\rm Num}}
\def\Irr{{\rm Irr}}
\def\Frac{{\rm Frac}}
\def\Sym{{\rm Sym}}
\def\TV{{\cal Z}}
\def\indlim{\underrightarrow{\lim}}
\def\prolim{\underleftarrow{\lim}}
\def\LE{{\rm LE}}
\def\LEG{{\rm LEG}}
\def\MFT{{\mathfrak T}}
\def\Exc{{\rm Exc}}
\def\Crit{{\rm Crit}}
\def\MFW{{\mathfrak W}}
\def\Transp{{\rm Transp}}
\def\min{{\rm min}}
\def\DB{\mathfrak{DB}}
\def\max{{\rm max}}


\def\RHom{{\rm RHom}}
\def\rRHom{{\mathcal RHom}}
\def\rHom{{\mathcal Hom}}
\def\dotimes{{\overline{\otimes}}} 
\def\Ext{{\rm Ext}}
\def\rExt{{\mathcal Ext}}
\def\Tor{{\rm Tor}}
\def\rTor{{\mathcal Tor}}
\def\SP{{\mathfrak S}}
\def\WR{{\mathfrak R}}

\def\H{{\rm H}}
\def\D{{\rm D}}
\def\Del{{\mathfrak D}}
\def\alg{{\rm alg}}
\def\sh{{\rm sh}}
 

 \newtheorem{theor}{Theorem}[section]
 \newtheorem{prop}[theor]{Proposition}
 \newtheorem{propdef}[theor]{Proposition-Definition}
 \newtheorem{sublemma}[theor]{sublemma}
 \newtheorem{cor}[theor]{Corollary}
 \newtheorem{lemma}[theor]{Lemma}
  \newtheorem{lemme}[theor]{Lemme}
 \newtheorem{sublem}[theor]{sub-lemma}
 \newtheorem{defin}[theor]{Definition}
 \newtheorem{conj}[theor]{Conjecture}
 \newtheorem{quest}[theor]{Question}
  \newtheorem{remark}[theor]{Remark}
    \newtheorem{lemma-def}[theor]{Lemma-Definition}

 
 \parindent=0pt

 \author{Damian R\"OSSLER\footnote{Mathematical Institute, 
University of Oxford, 
Andrew Wiles Building, 
Radcliffe Observatory Quarter, 
Woodstock Road, 
Oxford
OX2 6GG, 
United Kingdom}}
 \title{Strongly semistable sheaves and the Mordell-Lang conjecture over function fields}
\maketitle


\begin{abstract}
We give a new proof of the Mordell-Lang conjecture in positive characteristic, in the situation 
where the variety under scrutiny is a smooth subvariety of an abelian variety. Our proof is based on the theory of semistable sheaves in positive characteristic, in particular on  Langer's theorem that the Harder-Narasimhan filtration of sheaves becomes 
strongly semistable after a finite number of iterations of Frobenius pull-backs. The interest of this proof is that it provides 
simple effective bounds (depending on the degree of the canonical line bundle) for the degree of the isotrivial finite cover whose existence is predicted by the Mordell-Lang conjecture. We also present a conjecture on the Harder-Narasimhan filtration of the cotangent bundle of 
a smooth projective variety of general type in positive characteristic and a conjectural refinement of the Bombieri-Lang conjecture in positive characteristic. 
\end{abstract}

\section{Introduction}

Let $B$ be an abelian variety over an algebraically closed field $F$ of 
characteristic $p>0$. Let $Y$ be an integral closed subscheme of $B$. 
Let $\Lambda\subseteq B(F)$ be a subgroup. Suppose that $\Lambda\otimes_\mZ\mZ_{(p)}$ is 
a finitely generated $\mZ_{(p)}$-module (here, as is customary, we write $\mZ_{(p)}$ for the localization 
of $\mZ$ at the prime $p$). Suppose that $\Stab(Y)=0$. Here $\Stab(Y)=\Stab_B(Y)$ is the translation stabilizer of $Y$. 
This is the closed subgroup scheme  of $B$, which is characterized uniquely 
by the fact that for any scheme $S$ and any morphism $b: S\to B$, translation by $b$ on the
product $B\times_F S$ maps the subscheme $Y\times S$ to itself if and
only if $b$ factors through $\Stab_B(Y)$.  Its existence is proven in 
\cite[exp.~VIII, Ex.~6.5~(e)]{SGA3-2}.

The Mordell-Lang conjecture for $Y$ and $B$ is the following statement. 

\begin{theor}[Mordell-Lang conjecture for abelian varieties; Hrushovski \cite{Hrushovski-Mordell-Lang}]
Suppose that $Y\cap\Lambda$ is Zariski dense in $Y$. Then there is a projective variety $Y'$ over a finite subfield $\mF_{p^r} \subseteq F$ and a finite and surjective morphism $h: Y'_F\to Y$.

\label{MLtheor}
\end{theor}

Theorem \ref{MLtheor} was first proven by Hrushovski in \cite{Hrushovski-Mordell-Lang} 
using model-theoretic methods and other proofs were given in \cite{Rossler-MMML}, \cite{Ziegler-Mordell-Lang} and \cite{Benoist-Bouscaren-Pillay-MLMM}.

The main results of the present article are the following results, which we 
shall prove simultaneously.

\begin{theor} Suppose that $Y\cap\Lambda$ is Zariski dense in $Y$. 
Suppose that $Y$ is smooth over $F$.\\ If
$$
p>\dim({Y})^2\int_Y{\rm c}^1(\Omega_{Y})^{\dim(Y)}
$$
then there is a projective variety $Y'$ defined over a finite subfield $\mF_{p^r} \subseteq F$ and an 
isomorphism \mbox{$h: Y'_F\simeq Y$.}
\label{thflash}
\end{theor}

\begin{remark}\rm A weaker (but more cumbersome) inequality than $
p>\dim({Y})^2\int_Y{\rm c}^1(\Omega_{Y})^{\dim(Y)}
$, which also implies the conclusion of Theorem \ref{thflash} is given in \refeq{eq8} below.
\end{remark}

\begin{theor}
Suppose that $Y\cap\Lambda$ is Zariski dense in $Y$. 
Suppose that $Y$ is smooth over $F$. If $\Omega_Y$ is strongly 
semistable with respect to $\det(\Omega_Y)$ then there is a projective variety $Y'$ defined over a finite subfield $\mF_{p^r} \subseteq F$ and an 
isomorphism\,\,\,\, $h: Y'_F\simeq Y$.
\label{thflashp}
\end{theor}
Note that $\det(\Omega_Y)$ is ample by \cite[Lemma 6]{Abra-Sub}. A vector bundle $V$ on $Y$ is strongly semistable 
if $F_Y^{\ast,\circ r}V$ is semistable for all sufficiently large $r$. Here 
$F_Y$ is the absolute Frobenius endomorphism of $Y$. See section \ref{VBPC} below for more details.

In other words, if $Y$ is smooth and either $p>\dim({Y})^2\int_Y{\rm c}^1(\Omega_{Y})^{\dim(Y)}
$ or $\Omega_Y$ is strongly semistable with respect to $\Omega_Y$, then $h$ may be taken to be an isomorphism in Theorem \ref{MLtheor}. In particular, $h$ may be taken to be an isomorphism 
in Theorem \ref{MLtheor} if $\dim(Y)=1$. This is an early result of Samuel (see \cite{Samuel-Com}).

\begin{remark}\rm If $S$ is a smooth closed subvariety of an abelian variety
over a field of characteristic $0$ such that $\Stab(S)=0$ then 
$\Omega_S$ is semistable with respect to $\det(\Omega_S)$ by a classical result 
of Yau. Thus Theorem \ref{thflashp} may be viewed 
as an exact positive characteristic counterpart (when $Y$ is smooth) of the 
Mordell-Lang conjecture in characteristic $0$, because in the latter the analog 
of the morphism $h$ can always be taken to be an isomorphism. 
\end{remark}

\begin{remark}\rm 
Note the following interesting consequence of Theorem \ref{thflash}. 
Suppose that we are given a smooth projective subvariety with trivial stabiliser inside an abelian variety over the function field of a variety defined over a number field $L$. 
If this situation is reduced modulo a prime ideal $\mathfrak p$ of $\CO_K$, then for all but a finite number of such ideals, the Mordell-Lang conjecture \ref{MLtheor} holds for the reduced situation and $h$ can be taken to be an isomorphism (for any choice of $\Lambda$).
\end{remark}

Our proof of Theorems \ref{thflash} and \ref{thflashp} does not rely on existing proofs of Theorem \ref{MLtheor}. In particular, we provide in this text a complete proof of Theorem \ref{MLtheor} in the situation where $Y$ is smooth.
Our method of proof does not rely on the differential techniques of \cite{Hrushovski-Mordell-Lang} or 
on the Galois-theoretic techniques of \cite{Rossler-MMML}, \cite{Ziegler-Mordell-Lang} and \cite{Benoist-Bouscaren-Pillay-MLMM}. It is purely geometric and uses the theory of semistable sheaves on varieties of dimension $>1$, in particular on Langer's theorem that the Harder-Narasimhan filtration of sheaves becomes 
strongly semistable after a finite number of iterations of Frobenius pull-backs (see Theorem 
\ref{thlang} below). The possibility of giving a geometric condition on the prime number $p$ for the morphism $h$ to be an isomorphism is intrinsic to our method and the previous methods of proof of Theorem \ref{MLtheor} do not naturally lead to such bounds (be it only because the quantity 
$\int_Y{\rm c}^1(\Omega_{Y})^{\dim(Y)}$ never appears in them).

Here is a more detailed outline of the proof. We first show that if the assumptions of 
the Mordell-Lang conjecture are verified, then one can construct an infinite tower of 
torsors under vector bundles over $X$, which is trivialised by an infinite tower 
of finite surjective base-changes. This first step already appears in \cite{Rossler-MMML}. We then use the theory of semistable sheaves in positive characteristic, as developed by Langer in \cite{Langer-Semistable}, to show that this tower can be trivialised by a single 
finite purely inseparable morphism (see the proof of Theorem \ref{mainth} below). Here we need some simple facts about the 
slopes of the sheaf of differentials of a smooth subvariety of an abelian variety whose 
stabiliser is trivial (see Lemma \ref{ablem} below) and the key input from Langer's theory is Theorem \ref{thlang} below. We also 
need a cohomological result of Szpiro and Lewin-M\'en\'egaux (see Proposition \ref{propinj} below), which appears 
in their partial proof of the Kodaira vanishing theorem in positive characteristic. 
With this trivialisation in hand, we construct the variety  
$Y'$ whose existence is asserted in Theorem \ref{thflash}, by using Grothendieck's formal GAGA theorem. This method of descent is also used in \cite{Rossler-Rational}. The bound 
given in Theorem \ref{thflash} is deduced from a result of Langer (see Theorem 
\ref{laneff} below), which  compares 
the slopes of a torsion free sheaf with the slopes of its Frobenius pull-backs.

\medskip
Finally, we would like to state the following conjectures, which are suggested by 
our proof of Theorem \ref{thflash}.

\begin{conj}
Let $Z$ be a projective variety over $F$. 
 Suppose that $Z$ is smooth and that $\det(\Omega_Z)$ is an ample line bundle. Suppose also that $H^0(X,\Omega_Z^\vee)=0$. 
 Then $\bar\mu_{\min,\det(\Omega_Z)}(\Omega_Z)>0$.
 \label{conjam}
 \end{conj}
 Here $\bar\mu_{\min,\det(\Omega_Z)}(\cdot)$  refers to the Frobenius-stabilised 
minimal slope with respect to $\det(\Omega_Z).$ See section \ref{VBPC} below for the definition.

\begin{remark} \rm Lemma \ref{ablem} below shows that Conjecture \ref{conjam} is verified if $X$ can be embedded 
in an abelian variety. 
Also, note that it seems likely that there are "many" varieties $X$ satisfying the assumptions
of Conjecture \ref{conjam}, such that $\Omega_Z$ is strongly semistable with 
respect to $\det(\Omega_Z)$  (see the beginning of 
section \ref{VBPC} for this notion), which is a condition stronger than $\bar\mu_{\min,\det(\Omega_Z)}(\Omega_X)>0$). 
Indeed, recall that the cotangent bundle $\Omega_S$ of a smooth and projective variety $S$ over $\mC$ is semistable 
with respect to $\det(\Omega_S)$, if $\det(\Omega_S)$ is ample. This is a consequence of 
the main result of \cite{Yau-Calabi}. 
On the other hand, there is speculation (see for example \cite{Shepherd-Barron-SS} and the references therein) that in many situations the reduction modulo a prime number $p$ of a semistable sheaf is strongly semistable for "most" prime numbers $p$.
\end{remark}

Another conjecture concerns a possible generalisation of Theorem \ref{thflash} to a more 
general geometrical context.

\begin{conj} Let $Z$ be a smooth and projective variety over a finitely generated field $F_0$ of characteristic 
$p$. Suppose that $\det(\Omega_Z)$ is an ample line bundle. Suppose also that $Z(F_0)$ is Zariski dense in $Z$. 
If $$
p>\dim({Z})^2\int_Z{\rm c}^1(\Omega_{Z})^{\dim(Z)}
$$
then there is a projective variety $Z'$ over a finite subfield $\mF_{p^r} \subseteq \bar F_0$ and 
an isomorphism $$h: Z'_{\bar F_0}\simeq Z_{\bar F_0}.$$
\label{flashgen}
\end{conj}

Theorem \ref{thflash} and the Lang-N\'eron theorem show that Conjecture \ref{flashgen} holds if $Z$ can be embedded over $F_0$ 
into an abelian variety over $F_0$. In \cite{Rossler-Rational}, using some of the results of the present text, we show that Conjecture \ref{flashgen} holds if 
$\Omega_Z$ is ample and $F_0$ has transcendence degree one over its prime field. Conjecture 
\ref{flashgen} is an attempt to make the Bombieri-Lang over function fields 
in positive characteristic (see the introduction of \cite{Rossler-Rational} for a discussion) more precise when the variety under scrutiny is smooth. 

The structure of the text is the following. In section \ref{secmainr}, we shall formulate three more technical results, from which Theorem \ref{thflash} will be deduced. In particular, Corollary \ref{corSS} gives a strengthening of 
Theorem \ref{thflash} but is more complicated to formulate. The proofs of these results are given in 
section \ref{proofmainr}. In section \ref{VBPC}, we shall review the results from Langer's theory that 
we shall need, derive some simple consequences from it and also formulate the cohomological result of Szpiro and Lewin-M\'en\'egaux alluded to above. Finally in section \ref{lastsec}, the proof of 
Theorem \ref{thflash} is given.

The basic definitions for this article will be fixed in the next section. 

{\bf Basic notational conventions.} If $Z$ is a scheme of characteristic $p$, we write $F_Z:Z\to Z$ for the 
absolute Frobenius endomorphism of $Z$. The short-hand w.r.o.g. refers to 
"without restriction of generality".

\section{Main results}

\label{secmainr}

The definitions and notations given in this section will be used throughout the article. They differ 
from those used in the introduction, which will not be used again in the text.

 Let $k_0$ be an algebraically closed field of characteristic $p>0$ and let $U$ be a smooth variety  
over $k_0$. Let $\CA$ be an abelian scheme over $U$ and let $\CX\hookrightarrow\CA$ be a closed subscheme. We let $K_0$ be the function field of $U$ and 
let $A:=\CA_{K_0}$ (resp. $X:=\CX_{K_0}$) be the generic fibre of $\CA$ (resp. $\CX$). 

For all $n\geqslant 0$, we define 
$$\Crit^n(\CX,\CA):=[p^{n}]_*(J^{n}(\CA/U))\cap J^{n}(\CX/U).$$ 
Here $J^n(\bullet/U)$ refers to the $n$-th jet scheme of $\bullet$ over $U$. 
See \cite[par. 2]{Rossler-MMML} for this and some more explanations. The scheme $J^{n}(\CA/U)$ is naturally a commutative group scheme over $U$ and $[p^{n}]$ refers to the multiplication-by-$p^n$-morphism. The notation $[p^{n}]_*(J^{n}(\CA/U))$ refers to the scheme-theoretic image of 
$J^{n}(\CA/U)$ by $[p^{n}].$ There is a natural projective system of 
$U$-schemes
$$
\dots\to J^{n}(\CX/U)\stackrel{\Lambda^\CX_{n,n-1}}{\to}J^{n-1}(\CX/U)\to\dots \to J^{0}(\CX/U)=\CX.
$$
If $\CX$ is smooth over $U$ 
then the $J^{n-1}(\CX/U)$-scheme $J^{n}(\CX/U)$ carries a natural structure of torsor under 
the vector bundle $\Omega_{\CX/U}^\vee\otimes\Sym^n(\Omega_{U/k_0})$, where 
the vector bundles $\Omega_{\CX/U}^\vee$ and $\Sym^n(\Omega_{U/k_0})$ have been 
pulled back to $J^{n-1}(\CX/U)$ via the natural morphisms. 

In particular, we have projective system of 
$U$-schemes
$$
\dots\to\Crit^{2}(\CX,\CA)\to\Crit^1(\CX,\CA)\to\CX.
$$
and one can show that that the connecting morphisms in this system are finite. See \cite[par. 3.1]{Rossler-MMML} for this. We let $\Exc^n(\CA,\CX)\hookrightarrow\CX$ be the 
scheme-theoretic image of $\Crit^n(\CA,\CX)$ in $\CX$. We let $\Crit^n(A,X)$ (resp. $\Exc^n(A,X)\hookrightarrow X$) be the generic fibre of $\Crit^n(\CA,\CX)$ (resp. $\Exc^n(\CA,\CX)\hookrightarrow\CX$). 

Now fix once a for all an ample line bundle $M$ on $X_{\bar K_0}$. 
If $X$ is smooth over $K_0$ and $\Stab(X)=0$, a natural choice of 
an ample line bundle is $\det(\Omega_{X_{\bar K_0}})$. See \cite[Lemma 6]{Abra-Sub} for this. 

\begin{lemma-def}  Suppose that $X$ is smooth and geometrically connected over $K_0$ and 
that $\Stab(X)=0$. Then $\bar\mu_\min(\Omega_{X_{\bar K_0}})>0$ and 
$$
\mathfrak{DB}(X):=p^{\sup\{n\in\mN\ |\ H^0(X,F_X^{\ast,\circ n}\Omega^\vee_{X/K_0}\otimes\Omega_{X/K_0})\not=0\}}\leqslant{\bar\mu_{\rm max}(\Omega_{X_{\bar K_0}})\over\bar\mu_\min(\Omega_{X_{\bar K_0}})}\Big.
$$
\label{finlem}
\end{lemma-def}
Here again $\bar\mu_\min(\cdot)=\bar\mu_{\min,M}(\cdot)$  (resp. 
$\bar\mu_{\rm max}(\cdot)=\bar\mu_{{\rm max},M}(\cdot)$) refers to the Frobenius-stabilised 
minimal (resp. maximal) slope with respect to $M.$ See section \ref{VBPC} below for the definition.

Let now $\Gamma$ be a subgroup of $A(\bar K_0)$. Suppose that 
$$\Gamma=\Div^p(\Gamma_0):=\{\gamma\in A(\bar K_0)\ | \ \exists n\in\mN^*:\,(n,p)=1\,\,\&\,\,n\cdot\gamma\in\Gamma_0\}$$
where 
$\Gamma_0$ is a finitely generated subgroup of $A(K_0)$. In particular, $\Gamma\otimes\mZ_{(p)}$ is a finitely generated $\mZ_{(p)}$-module. 

\begin{theor}
 Suppose that $\CX$ is smooth over $U$ with geometrically connected fibres and 
suppose that $\Stab(X)=0$. Consider the statements: 
\begin{itemize}
\item[\rm (a)] For any $n\geqslant 0$ there is a $Q=Q(n)\in\Gamma_0$ such that $\Exc^n(A,X^{+Q})\hookrightarrow X$ is an isomorphism.
\item[\rm (b)] For any closed point $u_0\in U$, there is an $n_0=n_0(u_0)$ 
such that $p^{n_0}\leqslant\DB(X)$ and a finite and surjective morphism 
of $\widehat{\CO}_{u_0}$-schemes
$$
\iota=\iota_{u_0}:\CX_{u_0}^{p^{-n_0}}\times_{k_0}\widehat{\CO}_{u_0}\to
\CX_{\widehat{\CO}_{u_0}}
$$
of degree equal to $p^{\dim(X)n_0}$. 
\end{itemize}
Then {\rm (a)} implies {\rm (b).}
\label{mainth}
\end{theor}
Here ${U}_{u_0}$ is the spectrum of the local ring of $U$ at $u_0$ and $\widehat{U}_{u_0}$ is its completion. 
The notation $X^{+Q}$ refers to the pushforward by the addition-by-$Q$ morphism of the subscheme $X$ of $A$. The scheme $\CX_{u_0}$ is the $k_0$-scheme, which is the fibre of $\CX$ at $u_0$. The symbol $\CX_{u_0}^{p^{-r}}$ refers to the 
scheme obtained from $\CX_{u_0}$ by composing the structure map of $\CX_{u_0}$ 
with the $n$-th power $\Frob_{k_0}^{-1,\circ n}$ of the inverse of the 
absolute Frobenius morphism $\Frob_{k_0}$ of $\Spec\ k_0$ (recall that 
$\Frob_{k_0}$ is an automorphism because $k_0$ is perfect). 

Notice that the morphism $\iota$ must be flat by "miracle flatness" 
(see \cite[Th. 23.1]{Matsumura-Commutative}), since both source and target of $\iota$ are regular schemes.
By the degree of $\iota$, we mean as usual 
$$
{\rm deg}(\iota):=\rk(\iota_*(\CO_{\CX_{u_0}^{p^{-n_0}}\times_{k_0}\widehat{\CO}_{u_0}})),
$$
noting that $\iota_*(\CO_{\CX_{u_0}^{p^{-n_0}}\times_{k_0}\widehat{\CO}_{u_0}})$ is a locally free sheaf, since $\iota$ is flat. 

\begin{cor}
Suppose 
that $X_{\bar K_0}\cap \Gamma$ is dense in $X_{\bar K_0}$. Suppose also that 
$X$ is smooth and geometrically connected over $K_0$. 

Then there exists a smooth projective variety $X'$ over $k_0$ and a finite and surjective $K_0^\sep$-morphism $$h:X'_{K_0^\sep} \to (X/\Stab(X))_{K_0^\sep}$$  such that $$\deg(h)\leqslant\DB(X/\Stab(X))^{\dim(X/\Stab(X))}.$$

\label{corSS}
\end{cor}

\section{The geometry of vector bundles in positive characteristic}

\label{VBPC}

Let $L$ be an ample line bundle on a smooth and projective variety $Y$ over an algebraically closed field $l_0$. If $V$ is a torsion free coherent sheaf
on $Y$, we shall write $$\mu(V)=\mu_L(V)=\deg_L(V)/\rk(V)$$ for the slope 
of $V$ (with respect to $L$). Here $\rk(V)$ is the rank of $V$, which is 
the dimension the stalk of $V$ at the generic point of $Y$. Furthermore, 
$$\deg_L(V):=\int_X{\rm c}_1(V)\cdot{\rm c}_1(L)^{\dim(Y)-1}.$$
Here $c_1(\cdot)$ refers to the first Chern class with values in an arbitrary Weil cohomology 
theory and the integral sign $\int_X$ is a short-hand for the push-forward morphism to $\Spec\ l_0$ in that theory. 

Recall that $V$ is called semistable (with respect to $L$) if for every coherent subsheaf $W$ of $V$, 
we have $\mu(W)\leqslant\mu(V)$. The torsion free sheaf $V$ is called strongly semistable 
if $\char(l_0)>0$ and $F_Y^{\ast,\circ n}V$ is semistable for all $n\geqslant 0$.  

In general, there exists a filtration
$$
0=V_0\subseteq V_{1}\subseteq\dots\subseteq V_{r-1}\subseteq V_r=V
$$
of $V$ by subsheaves, such that the quotients $V_i/V_{i-1}$ are all semistable and such 
that the slopes $\mu(V_i/V_{i-1})$ are strictly decreasing for $i\geqslant 1.$ This filtration is unique 
and is called the Harder-Narasimhan (HN) filtration of $V$. We shall write
$$
\mu_\min(V):=\inf\{\mu(V_i/V_{i-1})\}_{i\geqslant 1}
$$
and
$$
\mu_\max(V):=\sup\{\mu(V_i/V_{i-1})\}_{i\geqslant 1}
$$
An important consequence of the definitions is the following fact: if $V$ and $W$ are two torsion 
free sheaves on $Y$ and $\mu_\min(V)>\mu_\max(W)$, then 
$\Hom_Y(V,W)=0$. 

For more on the theory of semistable sheaves, see the monograph \cite{Huybrechts-Lehn-The-geometry}. 

The following theorem will be a key input in our proof of Theorem \ref{mainth}. For the 
proof see \cite[Th. 2.7]{Langer-Semistable}.

\begin{theor}[Langer]
If $V$ is torsion free coherent sheaf
on $Y$ and $\char(l_0)>0$, then there exists $n_0\geqslant 0$ such that $F_Y^{\ast,\circ n}V$ has a strongly semistable 
HN filtration for all $n\geqslant n_0$. 
\label{thlang}
\end{theor}

If $V$ is a torsion free sheaf on $Y$ and $\char(l_0)>0$, we now define
$$
\bar\mu_\min(V):=\lim_{r\to\infty}\mu_\min(F_Y^{\ast,\circ r}V)/\char(l_0)^r
$$
and
$$
\bar\mu_\max(V):=\lim_{r\to\infty}\mu_\max(F_Y^{\ast,\circ r}V)/\char(l_0)^r.
$$
Note that Theorem \ref{thlang} implies that the sequences $\mu_\min(F_Y^{\ast,\circ r}V)/\char(l_0)^r$ 
(resp. $\mu_\max(F_Y^{\ast,\circ r}V)/\char(l_0)^r$) become constant when $r$ is sufficiently large, so 
the above definitions of $\bar\mu_\min$ and $\bar\mu_\max$ make sense. 

One can show that the sequence $\mu_\min(F_Y^{\ast,\circ r}V)$ 
(resp. the sequence \mbox{$\mu_\max(F_Y^{\ast,\circ r}V)$}) is decreasing 
(resp. increasing). 

We shall also need the following numerical estimate in the proof of Theorem \ref{thflash}. 
This is again a result of Langer, proved in \cite[Cor. 6.2]{Langer-Semistable} .

To formulate it, let 
$$
\alpha(V):=\max\{\mu_\min(V)-\bar\mu_\min(V),\bar\mu_\max(V)-\mu_\max(V)\}
$$
when $\char(l_0)>0$.

\begin{theor}[Langer]
If $\char(l_0)>0$, we have
$$
\alpha(V)\leqslant{\rk(V)-1\over\char(l_0)}\max\{\bar\mu_\max(\Omega_Y),0\}.
$$
\label{laneff}
\end{theor}

\begin{lemma-def}  Suppose that $\char(l_0)>0$. Suppose that 
$\bar\mu_\min(V)>0$. Then the quantity
$$
\mathfrak{DB}(V):=p^{\sup\{n\in\mN\ |\ H^0(X,F_Y^{\ast,\circ n}V^\vee\otimes \Omega_Y)\not=0\}}
$$
is finite and we have
$$
\mathfrak{DB}(V)\leqslant {\bar\mu_{\rm max}(\Omega_Y)\over\bar\mu_\min(V)}.
$$
\label{finlemprime}
\end{lemma-def}
\beginProof Notice that 
$$
H^0(X,F_Y^{\ast,\circ n}V^\vee\otimes \Omega_Y)\simeq
\Hom_X(F_Y^{\ast,\circ n}V, \Omega_Y)
$$
and furthermore, for any $r\geqslant 0$, there is a natural inclusion
$$
\Hom_X(F_Y^{\ast,\circ n}V, \Omega_Y)\subseteq 
\Hom_X(F_Y^{\ast,\circ(n+r)}V, F_Y^{\ast,\circ r}\Omega_Y)
$$
given by pulling back morphisms of vector bundles by $F_Y^{\ast,\circ r}$. 
Now by Theorem \ref{thlang}, we may choose $r$ sufficiently large so that $F_Y^{\ast,\circ r}V$ and $F_Y^{\ast,\circ r}\Omega_Y$ have 
Harder-Narasimhan filtrations with strongly semistable quotients.  Then we have
$$
\mu_\min(F_Y^{\ast,\circ(n+r)}V)=p^n\cdot\mu_\min(F_Y^{\circ r,*}V)
$$
and 
$$
\mu_\max(F_Y^{\ast,\circ(n+r)}\Omega_Y)=p^n\cdot\mu_\max(F_Y^{\circ r,*}\Omega_Y)
$$
Thus $\Hom_Y(F_Y^{\ast,\circ(n+r)}V, F_Y^{\circ r,*}\Omega_Y)=0$ 
if $$p^n\cdot\mu_\min(F_Y^{\ast,\circ r}V)>\mu_\max(F_Y^{\ast,\circ r}\Omega_Y).$$
Thus 
$$
\sup\{p^n \ |\,H^0(X,F_Y^{\ast,\circ n}V^\vee\otimes V)\not=0\}_{n\in\mN}\leqslant{\mu_\max(F_Y^{\ast,\circ r}\Omega_Y)\over\mu_\min(F_Y^{\ast,\circ r}V)}=
{\bar\mu_\max(\Omega_Y)\over \bar\mu_\min(V)}.\qed
$$

The next lemma is well-known but for lack of a bibliographical reference, we have included a proof.

\begin{lemma}
Let $V$ be a torsion free sheaf on $Y$. Suppose that $V$ is globally generated and of degree $0$ with respect to $L$. 
Then there exists an isomorphism $V\simeq \CO_Y^{\oplus\rk(V)}$. 
\end{lemma}
\beginProof
Let $\phi:\CO^{\oplus l}_Y\to V$ be a surjection, where $l$ is chosen as small as possible. 
Suppose that $\ker\ \phi\not=0$ (otherwise the Lemma is proven). 
Let $V_0=\ker\ \phi$. Then $\mu(V_0)=0$ and furthermore, since $\CO^{\oplus l}_Y$ is semistable, 
every semistable subsheaf of $V_0$ has slope $\leqslant 0$ and thus $V_0$ is also semistable. 
Now for any $i\in\{1,\dots,l\}$, let $\pi_i:V_0\to\CO_Y$ be the projection on the $i$-th coordinate. 
Choose $i_0\in\{1,\dots,l\}$ so that $\pi_{i_0}$ is non-vanishing. Then $\pi_{i_0}$ is surjective in codimension $2$, because 
otherwise, the degree of the image of $\pi_{i_0}$ would be $<0$, which would contradict the 
semistability of $V_0$. Now replace $V_0$ be a non-zero semistable subsheaf of $\ker \pi_{i_0}$ and repeat 
the above reasoning, unless $\pi_{i_0}$ is an isomorphism outside a closed subset of codimension at least $2$. Continuing in the same way, we  end up with a semistable torsion free sheaf  $M_0\subseteq \ker\ \phi\subseteq\CO^{\oplus l}$ of rank $1$, 
endowed with an arrow $M_0\to\CO_Y$, which is an isomorphism  outside a closed subset of codimension at least $2$. 
We thus obtain a complex $$\CO_Y|_{Y\backslash Y_0}\to \CO^{\oplus l}_Y|_{Y\backslash Y_0}\to 
V|_{Y\backslash Y_0},$$ where $Y_0$ is a closed subscheme of $Y$, which is of codimension at least $2$. 
Since $Y$ is normal, the arrow $\CO_Y|_{Y\backslash Y_0}\to \CO^{\oplus l}_Y|_{Y\backslash Y_0}$ extends uniquely to all of $Y$. We thus 
obtain a surjection $\CO^{\oplus l}_Y/\CO_Y\simeq \CO^{\oplus l-1}_Y\to V$. This contradicts 
the minimality of $l$ and proves the lemma. 
 \endProof

\begin{cor}
Let $V$ be a torsion free sheaf. Suppose that $V$ is globally generated. Then 
$V\simeq V_0\oplus\CO_Y^l$ for some $l\geqslant 0$ and for some torsion sheaf 
$V_0$ such that $\mu_\min(V_0)>0$. 
\label{cortriv}
\end{cor}
\beginProof Left as an exercise to the reader.\endProof
\begin{cor}
Let $V$ be a vector bundle over $Y$. 
Suppose that 

- for any surjective finite morphism $\phi:Y_0\to Y$, we have $H^0(Y_0,\phi^*V)=0;$

- $V^\vee$ is globally generated.

Then for any surjective finite morphism $\phi:Y_0\to Y$, such that $Y_0$ is smooth over $l_0$, we have $\mu_\min(\phi^*V^\vee)>0.$ In particular, if 
$\char(l_0)>0$ then $\bar\mu_\min(V^\vee)>0$. 
\label{cormmpos}
\end{cor}
\beginProof
The bundle $V^\vee$ is globally generated so  $\mu_\min(\phi^*V^\vee)\geqslant 0.$
Now to obtain a contradiction, suppose that $\phi^*V^\vee$ has a non-zero semistable quotient $Q$ of degree $0$. Then 
we have $\phi^*V^\vee\simeq Q_0\oplus\CO_{Y_0}^{\oplus l}$ for some $l>0$ 
by Corollary \ref{cortriv}. This implies that $\phi^*V$ has a non-vanishing section, which contradicts the assumptions.\endProof

The following elementary lemma is crucial to this article. The assumption that 
$Y$ is smooth over $l_0$ is not used in the next lemma. 

\begin{lemma}
Let 
$$
0\to V\to W\to N\to 0
$$
be an exact sequence of vector bundles on $Y$. 

Suppose that $W\simeq\CO_Y^{l}$ for some $l>0.$ 

Then $V^\vee$ is globally generated and for any dominant proper morphism $\phi:Y_0\to Y$, where $Y_0$ is integral, the morphism 
$$
\phi^*:H^0(Y,V)\to H^0(Y_0,\phi^*V)
$$
is an isomorphism.
\label{japanlem}
\end{lemma}
\beginProof
The fact that $V^\vee$ is globally generated follows from the fact that the natural dual 
map $W^\vee\to V$ is surjective. To prove the second statement, consider that 
we have a commutative diagram
\begin{diagram}
0 & \rTo & H^0(Y,V) & \rTo & H^0(Y,W) & \rTo & H^0(Y,N)\\
   &         &\dTo^{\phi^*} &     & \dTo^{\phi^*} &   & \dTo^{\phi^*}\\
0 & \rTo & H^0(Y_0,\phi^*V) & \rTo & H^0(Y_0,\phi^*W) & \rTo & H^0(Y_0,\phi^*N)  
\end{diagram}
In this diagram, all three vertical arrows are injective by construction. 
Furthermore, the middle vertical arrow is an isomorphism, also by construction. 
The five lemma now implies that the left vertical arrow is surjective.
\endProof

\begin{lemma}
Suppose that there is a closed $l_0$-immersion $i:Y\hookrightarrow B$, where $B$ is 
an abelian variety over $l_0$. Suppose that $\Stab_B(Y)=0$. 
Then $\Omega^\vee_Y$ is globally generated and for any dominant proper morphism \mbox{$\phi:Y_0\to Y$,} where $Y_0$ is integral,  we have $H^0(Y_0,\phi^*\Omega_Y^\vee)=0$. Furthermore, 
we have $\mu_\min(\Omega_Y)>0$ and if $\char(l_0)>0$, we have $\bar\mu_\min(\Omega_Y)>0$.
\label{ablem}
\end{lemma}
\beginProof 
We have an exact sequence 
$$
0\to\Omega_Y^\vee\to i^*\Omega_B^\vee\to N_{Y/B}\to 0
$$
where $N_{Y/B}$ is the normal bundle of $Y$ in $B$. 
Furthermore, since $\Stab_B(Y)=0$, we have $H^0(Y,\Omega_Y^\vee)=0$. 
Remembering that $\Omega_B$ is a trivial bundle, the lemma now follows from Lemma \ref{japanlem} and Corollary \ref{cormmpos}.\endProof

In the following lemma, the smoothness assumption on $Y$ is not used either. 
The proof of the following lemma is extracted from \cite[p. 49, before Prop. 3]{Martin-Deschamps-Proprietes}, where the argument is attributed to 
Moret-Bailly. 

\begin{lemma}
Suppose given a vector bundle $V$ on $Y$ with the following property: 
if $\phi:Y_0\to Y$ is a surjective and finite morphism and $Y_0$ is integral, then we have 
$H^0(Y_0,\phi^*V)=0$. 

Let $f:T\to Y$ be a torsor under $V$ and let $Z\hookrightarrow T$ be a closed immersion. 
Suppose that $f|_Z:Z\to Y$ is finite and surjective and that $Z$ is integral. 

Then $f|_Z$ is generically purely inseparable.
\label{radlem}
\end{lemma}
\beginProof  
Let $f:T\times_Y T\to Y$. 
We consider the scheme $T\times_Y(T\times_Y T)$. Via the projection 
on the second factor $T\times_Y T$, this scheme is naturally a torsor under the vector bundle $f^*V$ . This torsor has two sections: 

- the section $\sigma_1$ defined by the formula $t_1\times t_2\mapsto t_1\times (t_1\times t_2)$;

- the section $\sigma_2$ defined by the formula $t_1\times t_2\mapsto t_2\times (t_1\times t_2)$.

Since $T\times_Y(T\times_Y T)$ is a torsor under $f^*V$, there is a section 
$s\in H^0(T\times_Y T,f^*V)$ such that $\sigma_1+s=\sigma_2$ and by construction 
$s(t_1\times t_2)=0$ iff $t_1=t_2$. In other words, $s$ vanishes precisely on the 
the diagonal of $T\times_Y T$. 

Consider now the closed immersion $Z\times_Y Z\hookrightarrow T\times_Y T$. 
Suppose to obtain a contradiction that $f|_Z$ is not generically 
purely inseparable. Then there is an irreducible component $C$ of $Z\times_Y Z$, which 
is not contained in the diagonal and such that $f|_C:C\to Y$ is dominant and hence surjective. 

Indeed, if $f|_Z$ is not generically 
purely inseparable, then there is by constructibility an open subset $U\subseteq Y$, such 
that for any closed point $u\in U$, there is a point $P(u)\in Z_u\times_u Z_u$ such that 
$P(u)$ is not contained in the diagonal of $Z_u\times_u Z_u$. Hence there is an irreducible 
component of $Z\times_Y Z$, which does not coincide with the diagonal and furthermore there is one, which dominates $U$ for otherwise not every $P(u)$ would be contained in 
an irreducible component of $Z\times_Y Z$. 

Now consider $f|_C^*V$. By construction the section $s|_C\in H^0(C,f|_C^*V)$ does not 
vanish. This contradicts the assumption on $V$.
\endProof

We now quote a result proved in \cite[exp. 2, Prop. 1]{Szpiro-Seminaire-Pinceaux}.

\begin{prop}[Lewin-M\'en\'egaux, Szpiro] Suppose that $\char(l_0)>0$. 
If \mbox{$H^0(Y,F^*_Y(V)\otimes\Omega_Y)=0$} then the natural map of abelian groups 
$$H^1(Y,V)\to H^1(Y,F_Y^*V)$$ is injective.
\label{propinj}
\end{prop}

\begin{cor}Suppose that $\char(l_0)>0$.  Let $V$ be a vector bundle over $Y$. 
Suppose that 

- for any surjective finite morphism $\phi:Y_0\to Y$, where $Y_0$ is integral, we have $H^0(Y_0,\phi^*V)=0;$

- $V^\vee$ is globally generated. 

Then there is an $n_0\in\mN$ such that 
$H^0(S,F^{n,*}_Y(V)\otimes \Omega_Y)=0$ for all $n>n_0$ and we might choose 
$n_0\leqslant \DB(V)$. 

Furthermore, let $T\to Y$ be a torsor under $F^{n_0,*}_Y(V)$. 
Let $\phi:Y'\to Y$ be a finite surjective morphism and suppose that $Y'$ is integral. Then the map 
$$H^1(Y,F^{n_0,*}_Y(V))\to H^1(Y',\phi^*(F^{n_0,*}_Y(V)))$$
is injective.
\label{cortor}
\end{cor}
\beginProof\, (of Corollary \ref{cortor}). 
The existence of $n_0$ is a consequence of Corollary \ref{cormmpos} and Lemma 
\ref{finlemprime}. The upper bound for $n_0$ is also a consequence of Lemma \ref{finlemprime}. 

For the second assertion, notice that 
by Lemma \ref{radlem}, we may assume w.r.o.g. that $\phi$ is generically purely inseparable. 
Let $H$ be the function field of $Y$ and let $H'|H$ be the (purely inseparable) function field extension 
given by $\phi$. Let $k_0>0$ be sufficiently large so that the extension $H'|H$ factors through 
the extension $H^{p^{-k_0}}|H$. We may suppose w.r.o.g. that $Y'$ is a normal scheme, since we may replace $Y'$ by its normalization without restriction of generality. 
On the other hand the morphism $F_Y^{k_0}:Y\to Y$ gives a presentation of $Y$ as its own 
normalization in  $H^{p^{-k_0}}.$ Thus there is a natural factorization $Y\to Y'\stackrel{\phi}{\to} Y$, where the composition of the two arrows is given by  $F_Y^{k_0}$. Now by Proposition \ref{propinj} 
 there is a natural 
injection $H^1(Y,F^{k_0,*}_Y(V))\hookrightarrow H^1(Y,F_Y^{k_0,*}(F^{n_0,*}_Y(V)))$. 
Thus there is an injection $H^1(Y,F^{n_0,*}_Y(V))\to H^1(Y',\phi^*(F^{n_0,*}_Y(V))).$
\endProof

\section{Proof of  Lemma \ref{finlem}, Theorem \ref{mainth} and Corollary \ref{corSS}}

 \label{proofmainr}

{\bf Proof of Lemma \ref{finlem}.} Follows from Lemma \ref{finlemprime} 
and Lemma \ref{ablem}. 

\medskip
{\bf Proof of Theorem \ref{mainth}.} Let $Q\in\CA(U)$. Consider the infinite commutative diagram 
of $\CX$-schemes
\begin{diagram}
\dots & \rTo & \Crit^{2}(\CX^{+Q},\CA)& \rTo & \Crit^1(\CX^{+Q},\CA) &\rTo & \CX\\
  &       & \dInto                   &        & \dInto                  &       & \dTo^{=}\\
  \dots & \rTo & J^{2}(\CX/U)& \rTo & J^1(\CX/U) &\rTo & \CX\\
\end{diagram}
For any $n\geqslant 0$, we shall write 
\begin{diagram}
\dots & \rTo & \Crit^{2}(\CX^{+Q},\CA)^{(p^n)}& \rTo & \Crit^1(\CX^{+Q},\CA)^{(p^n)} &\rTo & \CX\\
  &       & \dInto                   &        & \dInto                  &       & \dTo^{=}\label{diagfund}\\
  \dots & \rTo & J^{2}(\CX/U)^{(p^n)}& \rTo & J^1(\CX/U)^{(p^n)} &\rTo & \CX\\
\end{diagram}
for the diagram obtained by pulling back the original diagram by $F_\CX^{\ast,\circ n}$.
Let 
$$
n_0:=\sup\{n\in\mN^*\ |\ H^0(X,F_X^{\ast,\circ n}\Omega^\vee_{X/K_0}\otimes\Omega_{X/K_0})\not=0\}.
$$
Suppose that (a) in Theorem \ref{mainth} is satisfied. We shall 
study diagram \refeq{diagfund} in the case where $n=n_0$. 
Now fix any $m>1$ and choose some $Q=Q(m)\in\CA(U)$ such that $\Exc^{m}(A,X^{+Q})\hookrightarrow X$ is an isomorphism. 
This is possible by assumption. By construction, 
the morphism $$\Crit^{m}(\CX^{+Q},\CA)^{(p^{n_0})}\to\CX$$ is then surjective. Choose 
an irreducible component $\Crit^{m}(\CX^{+Q},\CA)^{(p^{n_0})}_0\hookrightarrow 
\Crit^{m}(\CX^{+Q},\CA)^{(p^{n_0})}$, which dominates $\CX$. Endow 
 $\Crit^{m}(\CX^{+Q},\CA)^{(p^{n_0})}_0$ with its induced reduced scheme structure and for 
any $l<m$, let $\Crit^{l}(\CX^{+Q},\CA)^{(p^{n_0})}_0\hookrightarrow 
\Crit^{l}(\CX^{+Q},\CA)^{(p^{n_0})}$ be the irreducible component obtained 
by direct image from $\Crit^{m}(\CX^{+Q},\CA)^{(p^{n_0})}_0$. 

Now notice that by Corollary \ref{cortor}  and Lemma \ref{ablem} the base-change of the 
$F_X^{\ast,\circ n_0}(\Omega^\vee_{X/K_0}\otimes\Omega_{K_0/k_0})$-torsor $J^1(X/K_0)^{(p^{n_0})}\to X$ to $\bar K_0$ is trivial  and it is thus a trivial 
torsor. 
Let \mbox{$\sigma:X\to J^1(X/K_0)^{(p^{n_0})}$} be a section. The datum of the composed morphism 
$$\Crit^{1}(X^{+Q},A)^{(p^{n_0})}_0\to X\stackrel{\sigma}{\to}J^1(X/K_0)^{(p^{n_0})}$$ is equivalent 
to the datum of a section of the pull-back of $\Omega_{X/k_0}^\vee$ to 
$\Crit^{1}(X^{+Q},A)^{(p^{n_0})}_0$, which must vanish by Lemma \ref{ablem} (note that $\Omega_{K_0/k_0}$ is a trivial bundle). Hence 
the morphism  $\Crit^{1}(X^{+Q},A)^{(p^{n_0})}_0\to X$ is an isomorphism and 
is the image of $\sigma$. In particular, if $J^1(X/K_0)^{(p^{n_0})}\to X$ has a section over $X$, this section is unique. Furthermore, by Zariski's main theorem, the morphism $\Crit^{1}(\CX^{+Q},\CA)^{(p^{n_0})}_0\to \CX$ is 
an isomorphism. We now repeat this reasoning for the restriction to $\Crit^{1}(\CX^{+Q},\CA)^{(p^{n_0})}_0$ of the 
$F_X^{\ast,\circ n_0}(\Omega^\vee_{X/K_0}\otimes\Sym^2(\Omega_{K_0/k_0}))$-torsor $J^2(X/K_0)^{(p^{n_0})}\to J^1(X/K_0)^{(p^{n_0})}$ and 
we conclude that $$\Crit^{2}(\CX^{+Q},\CA)^{(p^{n_0})}_0\to\Crit^{1}(\CX^{+Q},\CA)^{(p^{n_0})}_0$$ is an 
isomorphism. Continuing this way, we see that in the whole tower
$$
\Crit^{m}(\CX^{+Q},\CA)^{(p^{n_0})}_0\to\Crit^{m-1}(\CX^{+Q},\CA)^{(p^{n_0})}_0\to\dots\to 
\Crit^{1}(\CX^{+Q},\CA)^{(p^{n_0})}_0\to\Crit^{1}(\CX^{+Q},\CA)^{(p^{n_0})}_0\to\CX
$$
the connecting morphisms are all isomorphisms. Letting $m\to\infty$, we obtain 
an infinite commutative diagram (\ref{infch}):
\begin{diagram}
\dots&\rTo&\CX^{m+1}&\rTo&\CX^{m}&\rTo&\dots&\rTo& 
\CX^1&\rTo&\CX\\
&&\dInto&&\dInto&&&& 
\dInto&&\dTo^{=}\\
\dots&\rTo&J^{m+1}(\CX/U)&\rTo&J^{m}(\CX/U)&\rTo&\dots&\rTo& 
J^{1}(\CX/U)&\rTo&\CX
\label{infch}
\end{diagram}
where all the morphisms $\CX^m\to \CX$ are isomorphisms.

Now choose a closed point $u_0\in U$.  View $u_0$ as a closed subscheme of $U$.  For any $i\geqslant 0$, let $u_i$ be the $i$-th infinitesimal neighborhood of $u_0\simeq \Spec\, k_0$ in $U$ (so that there is no ambiguity of notation for $u_0$). Notice that $u_i$ has a natural structure of $k_0$-scheme. Recall that by the definition of the jet scheme (see \cite[sec. 2]{Rossler-MMML}), 
the scheme $J^m(\CX/U)_{u_0}$ represents the functor 
on $k_0$-schemes
$$
T\mapsto{\rm Mor}_{u_m}(T\times_{k_0}{u_m},\CX_{u_m}).
$$
Thus the infinite chain \refeq{infch} gives rise to morphisms 
\begin{equation}
\CX_{u_0}^{(p^{-n_0})}\times_{k_0} u_m\to \CX_{u_m}
\label{formmor}
\end{equation}
compatible with each other under base-change. In particular, base-change to $u_0$ gives 
$F^{n_0}_{\CX_{u_0}}$. 

View the $\widehat{U}_{u_0}$-schemes $\CX_{u_0}^{(p^{-n_0})}\times_k \widehat{U}_{u_0}$ and $\CX_{\widehat{U}_{u_0}}$ as formal 
schemes over $\widehat{U}_{u_0}$ in the next sentence. The family
of morphisms \refeq{formmor} provides us with a morphism of formal schemes 
$$
\CX_{u_0}^{(p^{-n_0})}\times_k \widehat{U}_{u_0}\to\CX_{\widehat{U}_{u_0}}
$$
and since both schemes are projective over $\widehat{U}_{u_0}$, Grothendieck's GAGA theorem shows that this morphism of 
formal schemes comes from  a unique morphism of 
schemes
$$
\iota:\CX_{u_0}^{(p^{-n_0})}\times_k \widehat{U}_{u_0}\to\CX_{\widehat{U}_{u_0}}.
$$
By construction the morphism $\iota$ specializes to $F^{n_0}_{\CX_{u_0}}$ at the closed point $u_0$ of $\widehat{U}_{u_0}$. Since $F^{n_0}_{\CX_{u_0}}$ has finite fibres and $\iota$ is proper (since both source and target are projective over $\widehat{U}_{u_0}$), the morphism 
$\iota$ is quasi-finite by semicontinuity of fibre dimension and it is thus finite by Zariski's main theorem. The morphism 
$\iota$ is also flat by "miracle flatness", since both source and target are regular (see \cite[Th. 23.1]{Matsumura-Commutative}). 

Thus
$${\rm deg}(\iota)=p^{\dim(X)n_0}.$$
 Finally, $p^{n_0}\leqslant\DB(X)$ by  Lemma \ref{finlem}.

\medskip
{\bf Proof of Corollary \ref{corSS}.} We may replace $X$ by $X/\Stab(X)$ without restriction of generality in the statement 
of Corollary \ref{corSS}. Thus we may (and do) assume that $\Stab(X)=0$. 
Notice that by construction, for any $n\geqslant 1$, the natural homomorphism of groups
$$
\Gamma_0/p^n\Gamma_0\to \Gamma/p^n\Gamma
$$
is a surjection. Furthermore, $\Gamma_0/p^n\Gamma_0$ is finite since 
$\Gamma_0$ is finitely generated. Hence, using the assumptions of 
Corollary \ref{corSS}, we see that for any $n\geqslant 1$, there exists 
$Q=Q(n)\in\Gamma_0$, such that $X^{+Q(n)}\cap p^n\Gamma$ is dense 
in $X^{+Q}$. This implies that $\Exc^n(A,X^{+Q(n)})\hookrightarrow X$ 
is an isomorphism (see \cite[par. 3.2]{Rossler-MMML} for more details or this). 
Now applying Theorem \ref{mainth} (b), we obtain 
a surjective and finite morphism of $\widehat{\CO}_{u_0}$-schemes
$$
\CX_{u_0}^{p^{-n_0}}\times_{k_0}\widehat{\CO}_{u_0}
\to\CX_{\widehat{\CO}_{u_0}}
$$
for some closed point $u_0$ in $U$ (in fact any will do) and 
some $n_0\geqslant 0$ such that $p^{n_0}\leqslant\DB(X)$. Let $\widehat{K}_0$ be the fraction field of  ${\widehat{\CO}_{u_0}}$.  

Since $k_0$ is an excellent field, we know that the field extension $\widehat{K}_0|K_0$ is separable. On the other hand 
the just constructed finite and surjective morphism $\CX_{u_0}\times_{k_0}\widehat{K}_0\to X_{\widehat{K}_0}$ is defined 
over a finitely generated (as a field over $K$) subfield $K'_0$  of $\widehat{K}$. The field extension $K'_0|K$ is then still separable (because the extension $\widehat{K}_0|K_0$ is separable) and thus by the theorem on separating transcendence bases, there exists a variety $U'/K_0$, which is smooth over $K_0$ and whose function field is $K'_0$. Furthermore, possibly replacing $U'$ by one of its open subschemas, we may assume that 
the morphism $\CX_{u_0}\times_k {K}'_0\to X_{{K}'_0}$ extends to a finite and surjective morphism 
$$\alpha:\CX_{u_0}\times_{k_0} U'\to X_{U'}.$$ 
Let $P\in U'(K^\sep_0)$ be a $K^\sep_0$-point over $K$ (the set $U'(K^\sep_0)$ is not empty because 
$U'$ is smooth over $K_0$). The morphism $\alpha_P$ is the morphism $h$ advertised in Theorem \ref{mainth} (b). The inequality 
$\deg(h)\leqslant\DB(X)^{\dim(X)}$ is verified by construction.

\section{Proof of Theorem \ref{thflash}}

\label{lastsec}

We shall use the shorthand $\Omega:=\Omega_Y$. We shall derive Theorem 
\ref{thflash} from the results of section \ref{secmainr}. We use the notation of that section. Let $k_0:=\bar\mF_p$ and let $K_0$ be a finitely 
generated extension of $k_0$, such that $Y$ admits a model over $K_0$. 
We define $X$ to be such a model and we choose a smooth variety 
$U$ over $k_0$, such that $\kappa(U)=K_0$ and such that 
$X$ extends to a smooth an projective scheme $\CX$ over $U$. 
We let $M:=\det(\Omega)$ (recall that $\det(\Omega)$ is ample). 
Corollary \ref{corSS} now implies Theorem \ref{thflash} (resp. Theorem \ref{thflashp}), provided we can 
show that $\mathfrak{DB}=1$ under the assumptions of Theorem \ref{thflash} (resp. Theorem \ref{thflashp}). If $\Omega$ is strongly semistable then  Lemma \ref{finlemprime} immediately implies that $\mathfrak{DB}=1$ so Theorem \ref{thflashp} is proven. In particular, to prove Theorem \ref{thflash} we may assume w.r.o.g. that $\Omega$ is not strongly semistable. 

We shall now use Langer's theorem \ref{laneff} to derive Theorem \ref{thflash}. 

Note first that by Lemma \ref{ablem}, we have $\bar\mu_\min(\Omega)>0$ and $\mu_\min(\Omega)>0$. In particular, we have $\deg(\Omega)\geqslant 1$. In view of Lemma \ref{finlemprime} again, it is now sufficient to show that the inequality
$$
p>d^2\deg(\Omega)
$$
ensures that the inequality
\begin{equation}
\bar\mu_\max(\Omega)<p\cdot\bar\mu_\min(\Omega).
\label{eq1}
\end{equation}
is verified. In view of Theorem \ref{laneff}, inequality \refeq{eq1} is implied by 
the inequality
\begin{equation}
\mu_\max(\Omega)+\alpha(\Omega)<p\cdot(\mu_\min(\Omega)-\alpha(\Omega))
\label{eq2}
\end{equation}
From the definitions, we have
\begin{equation}
\bar\mu_\max(\Omega)\leqslant \deg(\Omega)/2
\label{eq3}
\end{equation}
and 
\begin{equation}
\mu_\max(\Omega)\leqslant \deg(\Omega)
\label{eq4}
\end{equation}
(recall that $\bar\mu_\min(\Omega)\not=\bar\mu_\max(\Omega)$ since 
$\Omega$ is assumed not to be strongly semistable) and
\begin{equation}
\mu_\min(\Omega)\geqslant 1/d
\label{eq5}
\end{equation}
Thus by Theorem \ref{laneff}, inequality \refeq{eq2} is weaker than the inequality
\begin{equation}
\deg(\Omega)+{(1+p)(d-1)\over 2p}\deg(\Omega)<p/d
\label{eq6}
\end{equation}
which can be rewritten as
\begin{equation}
2p^2-(d^2+d)\deg(\Omega)p-(d^2-d)\deg(\Omega)>0
\label{eq7}
\end{equation}
The roots of the equation in $x$
\begin{equation}
2x^2-(d^2+d)\deg(\Omega)x-(d^2-d)\deg(\Omega)=0
\label{eqp}
\end{equation}
are 
$$
{1\over 4}\Big[(d^2+d)\deg(\Omega)\pm\sqrt{\deg(\Omega)^2(d^2+d)^2+8\deg(\Omega)(d^2-d)}\Big].
$$
One of these roots is $\leqslant 0$ and the other one is $\geqslant 0$. 
Thus the inequality \refeq{eq7} is equivalent to the inequality
\begin{equation}
p> {1\over 4}\Big[(d^2+d)\deg(\Omega)+\sqrt{\deg(\Omega)^2(d^2+d)^2+8\deg(\Omega)(d^2-d)}\Big].
\label{eq8}
\end{equation}
In particular, if 
$$
p>d^2\deg(\Omega)
$$
then inequality \refeq{eq8} is verified (to check this quickly, just notice that  
 the function 
$$
2x^2-(d^2+d)\deg(\Omega)x-(d^2-d)\deg(\Omega)
$$ evaluated 
at $x=d^2\deg(\Omega)$ is $\geqslant 0$). 


\end{document}